\documentclass{amsart}
\usepackage{amssymb, latexsym}
\usepackage{amscd}

\newtheorem{theorem}{Theorem}

\newtheorem{definition}{Definition}

\begin{document}

\title{A Proof of Simon's Conjecture}
\author{Rita Gitik}
\email{ritagtk@umich.edu}
\address{Department of Mathematics \\ University of Michigan \\ Ann Arbor, MI, 48109}  

\date{\today}

\begin{abstract}
We prove Simon's conjecture for $3$-manifolds. 
\end{abstract}

\maketitle


\section{Introduction}

A manifold $M$ is called a missing boundary manifold  if it can be 
embedded in a compact manifold $\overline{M}$  such that $\overline{M} - M$
is a closed subset of the boundary of $\overline{M}$. J. Simon gave the following definition in his venerable paper \cite{Si}.

\begin{definition}
If a $3$-manifold $M$ has the property that for any finitely generated 
subgroup $H$ of $\pi_1(M)$ the cover of $M$ corresponding to $H$ is a missing boundary manifold, then we shall say that $M$ has almost-compact coverings.
\end{definition}

Simon conjectured in \cite{Si} that every $P^2$-irreducible, compact, connected $3$-manifold has almost-compact coverings.

Thurston's Geometrization Conjecture, proved by Perelman, (\cite{B-B-B-M-P} and \cite{K-L}), implies that
every compact orientable $3$-manifold admits a decomposition
along a (possibly empty) collection of spheres, discs, incompressible annuli, and incompressible tori, such that each
of the resulting $3$-manifolds is either a Seifert manifold or hyperbolic.

Using modern terminology, Simon proved (Theorem 2.5 in \cite{Si}) that if $M$ is an orientable compact $3$-manifold satisfying Thurston's Geometrization Conjecture and all submanifolds of $M$ obtained from the decomposition of $M$ given by the Geometrization Conjecture have almost-compact coverings then $M$ has almost-compact coverings.

Our proof of Simon's conjecture follows from two theorems.

\begin{theorem}
Seifert manifolds have almost-compact coverings.
\end{theorem}

\begin{theorem}
Compact $3$-manifolds with hyperbolic interiors have almost-compact coverings.
\end{theorem}

\section{Proof of Theorem 1} 

\begin{definition} cf. \cite{Sc}.

An orientable Seifert manifold is a $3$-manifold $M$ foliated by disjoint circles, such that each circle has a fibered neighborhood in $M$ homeomorphic to a fibered solid torus.
\end{definition}

\textbf{Proof of Theorem 1.}

Any compact orientable Seifert manifold is finitely covered by a bundle over an orientable compact surface with fiber a circle.
Any cover of such a bundle is either a bundle over an orientable surface with fiber a circle or an orientable bundle over an orientable surface with fiber $\mathbf{R}$. If the fundamental group of the bundle is finitely generated, so is the fundamental group of the base space.
An orientable bundle over an orientable surface $F$ with fiber $\mathbf{R}$ is homeomorphic to $F \times \mathbf{R}$. As $\pi_1(F)$ is finitely generated, $F$  is a missing boundary surface with compactification which we denote by $\overline{F}$. So $F \times \mathbf{R}$ is a missing boundary manifold with compactification $\overline{F} \times \mathbf{I}$. A circle bundle over a missing boundary surface $F$ extends to a circle bundle over $\overline{F}$, which is compact. Therefore Seifert manifolds have almost-compact coverings, proving Theorem 1.

\section{Proof of Theorem 2} 

\begin{definition}
A hyperbolic $3$-manifold is the quotient of $\mathbf{H^3}$ by a discrete group of isometries acting freely.
\end{definition}

\begin{definition}
The convex hull of a hyperbolic $3$-manifold $M$ is the minimal convex submanifold $C(M)$ of $M$. It is characterized by the following property:
for any pair of points $x$ and $y$ in $C(M)$ any geodesic segment connecting $x$ and $y$ is contained in $C(M)$.
\end{definition}

\textbf{Proof of Theorem 2.}

If a compact $3$-manifold with hyperbolic interior has a boundary component which is not a torus, then it has infinite volume. The fundamental group of such a manifold has the finitely generated intersection property (f.g.i.p.), cf. \cite{He}, so the result follows from \cite{Si}, Theorem 3.7.

If all the boundary components of a compact $3$-manifold with hyperbolic interior are tori, then it has finite volume. In this case we use the following result, proved by W. Thurston in \cite{Thu} and R. Canary in \cite{Ca2}.

\begin{theorem} cf. \cite{Ca1}, Corollary 8.1.

If $M$ is a finite volume hyperbolic $3$-manifold and $H$ is a finitely generated subgroup of $\pi_1(M)$, then $H$ is either 
geometrically finite or a virtual fiber group.
\end{theorem}

If $H$ is a virtual fiber group, i.e. $M$ is finitely covered by a bundle over $\mathbf{S^1}$ with fiber $F$, and $\pi_1(F)=H$, then the cover of $M$ corresponding to $H$ is homeomorphic to $F \times \mathbf{R}$.
Hence that cover is a missing boundary manifold with compactification $\overline{F} \times \mathbf{I}$.

If $H$ is geometrically finite, i.e. an $\epsilon$-neighborhood of the convex core of the cover of $M$ corresponding to $H$ has finite volume for all positive $\epsilon$, then Theorem 4 implies that the cover of $M$ corresponding to $H$ is a missing boundary manifold, proving
Theorem 2.

\begin{theorem}
The cover of a finite volume hyperbolic $3$-manifold corresponding to a geometrically finite subgroup of its fundamental group is a missing boundary manifold.
\end{theorem}

\section{Proof of Theorem 4.}

Let $H$ be a subgroup of a group $G$ given by the  presentation $G= \langle X|R \rangle $. 
Let $K$ be the standard presentation $2$-complex of $G$, i.e. $K$ has one 
vertex, $K$ has an edge, which is a loop, for every generator $x \in X$, and $K$ has a $2$-cell
for every relator $r \in R$. The Cayley complex of $G$, denoted by $Cayley_2(G)$, is  
the universal cover of $K$. Denote by $Cayley_2(G,H)$ the cover of $K$ corresponding to
a subgroup $H$ of $G$.

\begin{definition} cf. \cite{Gi1}, \cite{Gi2}, \cite{Gi3}, and \cite{Mi}.

A finitely generated subgroup $H$ of a finitely presented group  
$G$ is tame in $G$ if for any finite subcomplex
$C$ of $Cayley_2(G,H)$ and for any component $C_0$ of $Cayley_2(G,H)-C$ the group
$\pi_1(C_0)$ is finitely generated.
\end{definition}

T. Tucker proved the following result in \cite{Tuc}.

\begin{theorem} cf. \cite{Tuc}, p.267.

Let $M_0$ be a compact orientable irreducible $3$-manifold, and $M$ be the cover of $M_0$
corresponding to  a finitely generated subgroup of $\pi_1(M_0)$. 
Then $M$ is a missing boundary manifold  if and only if $\pi_1(M)$ is tame in $\pi_1(M_0)$.
\end{theorem} 

M. Mihalik proved the following result in \cite{Mi}.

\begin{theorem} cf. \cite{Mi}, Theorem 4.

If $H$ is a quasiconvex subgroup of the automatic group $G$, then $H$ is tame in $G$.
\end{theorem}

We need the following result which was proved in \cite{Ep}, p. 266.

\begin{theorem} cf. \cite{Ep}, Theorem 11.4.1.

 The fundamental group of a geometrically finite hyperbolic $3$-manifold is biautomatic.
\end{theorem}

Benson Farb informed the author in a  private conversation that any geometrically finite subgroup of the fundamental group of a geometrically finite hyperbolic $3$-manifold is quasiconvex with respect to the biautomatic structure constructed in the proof of Theorem 5 in \cite{Ep}.

Therefore Theorem 4 follows from Theorem 6.

\bigskip

\textbf{Acknowledgment.}

The author would like to thank Benson Farb and Peter Scott for helpful discussions and Shmuel Weinberger for his support.

\end{document}